\documentclass[11 pt]{amsart}
\usepackage{amsmath}
\usepackage{amssymb}
\usepackage{amsthm}
\usepackage{dsfont}
\usepackage{fancyhdr}
\usepackage[all]{xy}
\usepackage{hyperref}
\usepackage{paralist}
\usepackage{tikz}
\usepackage{verbatim}
\usepackage[hmargin=1.2in,vmargin=1.2in]{geometry}

\newtheorem{theorem}{Theorem}[section]
\newtheorem{lemma}[theorem]{Lemma}
\theoremstyle{proposition}

\theoremstyle{corollary}

\theoremstyle{definition}

\numberwithin{equation}{section}

\newcommand{\cE}{{\mathcal E}}

\newcommand{\cF}{{\mathcal F}}
\newcommand{\cG}{{\mathcal G}}
\newcommand{\cH}{{\mathcal H}}

\newcommand{\Z}{\mathbb Z}

\newcommand{\cL}{{\mathcal L}}
\newcommand{\F}{\mathbb F}

\newcommand{\f}{\mathfrak f}

\newcommand{\g}{\mathfrak g}

\begin{document}
\title[]{On product decomposition}
\author[]{Ming-Deh A. Huang (USC, mdhuang@usc.edu)}
\address{Computer Science Department,University of Southern California, U.S.A.}
\email{mdhuang@usc.edu}

\urladdr{}
\date{\today}
\keywords{polynomial system, product decomposition, last fall degree}
\subjclass[2010]{13P10, 13P15}

\maketitle

\begin{abstract}
Given a finite set $W$ in $\bar{k}^n$ where $\bar{k}$ is the algebraic closure of a field $k$ one would like to  determine if $W$ can be decomposed as $\prod_{i=1}^n V_i$ where $V_i \subset \bar{k}$ under a linear transformation, that is, $W\stackrel{\lambda}{\to} \prod_{i=1}^n V_i$ where $\lambda\in Gl_n (\bar{k})$.  We assume that $W$ is presented as $W=Z(\mathcal{F})$, the zero set of a polynomial system $\mathcal{F}$ in $n$ variables over $k$.
We study algebraic characterization of such product decomposition.  For decomposition into component sets of the same cardinality we obtain a stronger characterization and show that the decomposition in this case is essentially unique (up to permutation and scalar multiplication of coordinates).  We investigate computational problems that arise from the decomposition problem.
\end{abstract}

\section{Introduction}
Suppose $W\subset\bar{k}^n$ is a finite set where $k$ is a field and $\bar{k}$ is the algebraic closure of $k$.
We say that $W$ has a product decomposition if
$W=\lambda(\prod_{i=1}^n V_i )$ where $V_i\subset\bar{k}$ is finite and $\lambda\in Gl_n (\bar{k})$.

We study algebraic characterizations of product decomposition.  In \S~\ref{alg-char} we begin with a simple characterization in Theorem~\ref{simple-main}.   We then provide a stronger characterization in Theorem~\ref{main} for the case where $W\subset\bar{k}^n$ has product decomposition
$W\stackrel{\lambda}{\to} \prod_{i=1}^n V_i$ with $\lambda\in Gl_n (\bar{k})$, $V_i\subset\bar{k}$ and $| V_i | = d$ for $i=1,\ldots, n$.  Moreover in Theorem~\ref{main-unique} we show that in this case product decomposition is unique up to permutation and scalar multiplication of coordinates, provided $d$ is not divisible by the characteristic of $k$.

Suppose
$W=Z(\cG)$, the zero set of a set $\cG$ consisting of $n$ polynomials $g_1,\ldots, g_n$ of degree bounded by $d$ in $\bar{k} [ x_1,\ldots,x_n]$.  Then by Theorem~\ref{main} $Z(\cG)$ has product decomposition into component sets of the same cardinality   if and only if
there exists $\mu\in Gl_n (\bar{k})$ such that letting $\cL$ be the linear vector space spanned by $\cG$, then the isomorphic linear vector space $\mu^* (\cL) = \{ h\circ\mu: h\in\cL\}$ has a basis
$f_i (x_i)$, $i=1,\ldots,n$, where each $f_i (x_i)$ is of degree $d$ without multiple roots.
Therefore such a decomposition exists for $W$ if and only if
\begin{equation}
\label{decompose}
\left( \begin{array}{c} g_1\\. \\. \\g_n\end{array} \right) =\rho\circ \left( \begin{array}{c} f_1\\. \\. \\f_n\end{array} \right)\circ\lambda
  \end{equation}
for some $f_i (x_i)$, $i=1,\ldots,n$, of degree $d$ without multiple roots, and $\lambda, \rho\in Gl_n (k)$.

Treating $f_i$, $\lambda,\rho$ as unknown, we arrive at a problem
of solving a system of polynomials of degree $d$ in $O(n^2)$ unknown and $O(n^d)$ conditions.  In general algebraic algorithms for solving polynomial systems take time exponential in the number of variables, which is $O(n^2)$ in this case (see for example \cite{LL}).
Whether the system characterized by Eq.~\ref{decompose} admits more efficient solution is an interesting open question.

In \S~\ref{compute} we study last fall degrees of polynomial systems associated with algebraic sets that admit product decomposition.  More generally we consider polynomial systems specified by a set
$\cG = \{ g_1,\ldots,g_n\}$ with $g_i\in k[x_1,\ldots,x_n]$ and $\deg (g_i) \le d$ for $i=1,\ldots,n$, such that
there exists $\mu\in Gl_n (k)$ where the linear space spanned by $g_i\circ\mu$, $i=1,\ldots,n$, contains some
$f_i \in k[x_i]$ with $\deg f_i \le d$ for $i=1,\ldots,n$.
We study the last fall degrees \cite{HKY,HKYY} of such polynomial
systems.
In Theorem~\ref{dG} we prove $d_{\cG} \le d$, where $d_{\cG}$ is the last fall degree of $\cG$.  If $d=O(1)$ then $d_{\cG} = O(1)$.
Finding $k$-rational solutions amounts to solving
the polynomial system $\cG\cup \cE$, where $\cE$ consists of the field equations $x_i^q - x_i$, $i=1,\ldots,n$ and $k = \F_q$.  In Theorem~\ref{dGE} we show that the last fall degree  $d_{\cG\cup\cE}\le d+q$.  In situations where the number of $k$-rational solutions is $O(1)$, it follows from Proposition 2.11 of \cite{HKYY} that the problem can be solved in time  $n^{O(d+q)}$, hence $n^{O(q)}$ if $d=O(1)$.
An improvement can be made in finding $k$-rational point by reducing the problem to solving a system $\cG\cup \cE'$ in $O(nm)$ variables and degree $O(dp)$, where $q=p^m$ and $\cE'$ consists of $x_{ij}^p-x_{i\ j+1}$ for $i=1,\ldots,n$ and $j=0,\ldots, m-1$, and $j+1$ is taken modulo $m$.  In Theorem~\ref{dGE} we show that the last fall degree  $d_{\cG\cup\cE'}\le dp$.  So if the number of $k$-rational solutions is $O(1)$, then it follows that the problem can be solved in time  $(nm)^{O(dp)}$, which is $(nm)^{O(1)}$ providing $d=O(1)$ and $p=O(1)$.

\section{Algebraic characterization}
\label{alg-char}
We begin with a simple algebraic characterization of product decomposition.
\begin{theorem}
\label{simple-main}
Suppose $W\subset\bar{k}^n$.
The following two statements are equivalent:
\begin{enumerate}
\item $W\stackrel{\lambda}{\to} \prod_{i=1}^n V_i$ with $\lambda\in Gl_n (\bar{k})$ and $V_i$ is a finite subset of $\bar{k}$ for $i=1,\ldots, n$.
\item
There exist polynomials $f_i (x_i)$, $i=1,\ldots,n$, each having distinct roots in $\bar{k}$, and $\lambda \in Gl_n (\bar{k})$ such that $W=Z(\{f_1\circ\lambda, \ldots,f_n\circ\lambda\})$.
\end{enumerate}

\end{theorem}

\ \\{\bf Proof} Suppose $W\stackrel{\lambda}{\to} \prod_{i=1}^n V_i$ with $\lambda\in Gl_n (\bar{k})$ and $V_i$ is a finite subset of $\bar{k}$ for $i=1,\ldots, n$.   Let $f_i (x_i ) = \prod_{\alpha\in V_i} (x_i -\alpha)$ for $i=1,\ldots,n$.
Then for $\alpha\in\bar{k}^n$, $\alpha\in W$ if and only if  $\lambda(\alpha)\in\prod_{i=1}^n V_i$, if and only if $f_i\circ\lambda (\alpha) = 0$ for $i=1,\ldots,n$.  So (1) implies (2).

Now assume (2).  Put $V = \prod_{i=1}^n V_i$ where $V_i$ is the set of roots of $f_i (x_i)$.  Then $\alpha\in W$ if and only if $f_i\circ \lambda(\alpha) = 0$ for $i=1,\ldots,n$ if and only if $\lambda(\alpha)\in V$.  So (2) implies (1).  $\Box$

Suppose $I$ is an ideal of $\bar{k}[x_1,\ldots,x_n]$. For $i\ge 0$, let $I_i$ denote the submodule of $I$ consisting of all $f\in I$ with $\deg f=i$.

For $S\subset \bar{k} [ x_1,\ldots, x_n]$ let $Z(S)$ denote the zero set of $S$, that is $Z(S)=\{ v\in \bar{k}^n: f(v)=0 \ \ \forall f\in S\}$.
For $V\subset\bar{k}^n$ let $I(V) = \{ f\in \bar{k}[x_1,\ldots,x_n]: f(v)=0 \ \ \forall v\in V\}$.

Theorem~\ref{main} and Theorem~\ref{main-unique} below provide a stronger characterization for the case of product decomposition into components of the same cardinality.  More precisely we restrict our attention to product decomposition of the form $W\stackrel{\lambda}{\to} \prod_{i=1}^n V_i$ with $\lambda\in Gl_n (\bar{k})$ and $V_i$ is a finite subset of $\bar{k}$ of the same cardinality $d$ for $i=1,\ldots, n$.

\begin{theorem}
\label{main}
Suppose $W\subset\bar{k}^n$ and $W=Z(\cG)$ where $\cG$ consists of $n$ polynomials $g_1,\ldots, g_n$ of degree bounded by $d$ in $\bar{k} [ x_1,\ldots,x_n]$.
Then the following are equivalent:
\begin{enumerate}
\item
$W\stackrel{\lambda}{\to} \prod_{i=1}^n V_i$ with $\lambda\in Gl_n (\bar{k})$, $V_i\subset\bar{k}$ and $| V_i | = d$ for $i=1,\ldots, n$.
\item
$g_1,\ldots,g_n$ are linearly independent and they form a basis of $I(W)_d$ as a $\bar{k}$-vector space.  Moreover there are $f_i (x_i)$, $i=1,\ldots,n$, of degree $d$ without multiple roots such that $f_i (x_i)$, $i=1,\ldots,n$, form a basis for the vector space spanned by $g_i\circ\mu$, $i=1,\ldots,n$, for some $\mu\in Gl_n (\bar{k})$.
\end{enumerate}
\end{theorem}

\begin{theorem}
\label{main-unique}
Suppose $W\subset\bar{k}^n$ and $| W | = d^n$.  Suppose $d$ is not divisible by the characteristic of $k$.  If $W\stackrel{\lambda}{\to} \prod_{i=1}^n V_i$ with $\lambda\in Gl_n (\bar{k})$, $V_i\subset\bar{k}$ and $| V_i | = d$ for $i=1,\ldots, n$, then such product decomposition is unique up to permutation and coordinate-wise scalar multiplication.
More explicitly
if $W\stackrel{\lambda'}{\to} \prod_{i=1}^n V'_i$ with $V'_i\subset\bar{k}$ and $| V'_i | = d$ for $i=1,\ldots, n$ then
$\lambda'=\rho\circ\lambda$ where $\rho (x_1,\ldots,x_n) = (\alpha_1 x_{\sigma_1},\ldots, \alpha_n x_{\sigma_n})$ for all $(x_1,\ldots,x_n)\in\bar{k}^n$, where $\alpha_1,\ldots,\alpha_n\in\bar{k}$ and $(\sigma_1,\ldots,\sigma_n)$ is a permutation of $(1,\ldots,n)$.
\end{theorem}

To prove the two theorems, we will need some technical preparation.

\begin{lemma}
\label{f=0}
Suppose $V_i$ is a finite set of $\bar{k}$ with $| V_i | = d_i$ for $i=1,\ldots,n$.
For $f(x_1,\ldots,x_n)\in\bar{k} [x_1,\ldots,x_n]$ with $\deg_{x_i} f < d_i$ for $i=1,\ldots,n$, if $f(\alpha_1,\ldots,\alpha_n) = 0$ for all $(\alpha_1,\ldots,\alpha_n)\in\prod_{i=1}^n V_i$, then $f=0$.
\end{lemma}
\ \\{\bf Proof} The lemma can be proved by induction on $n$.  The statement is true for $n=1$ since a polynomial in $x_1$ of degree less than $d_1$ cannot have $d_1$ roots unless it is 0.

\ \\Suppose the statement is true for $n-1$.

\ \\Write $f = \sum_{i=0}^{d_1 - 1} f_i x_1^i$ where $f_i\in\bar{k}[x_2,\ldots,x_n]$.
For $j > 1$ since $\deg_{x_j} f < d_j$ and $\deg_{x_j} f = \max_{i=0,\ldots,d_1 -1} \deg_{x_j} f_i$, it follows that
$\deg_{x_j} f_i < d_j$ for $i=0,\ldots,d_1 -1$.

\ \\Let $(\beta_2,\ldots,\beta_n)\in\prod_{i=2}^n V_i$.  Then by assumption $f (x_1,\beta_2,\ldots,\beta_n)$ vanishes at
every $\beta_1\in V_1$, hence it has at least $d_1$ roots.  Since its degree in $x_1$ is no greater than $\deg_{x_1} f$, which is less than $d_1$, we conclude that it is 0.  So $f_i (\beta_2,\ldots,\beta_n)=0$ for $i=0,\ldots,d_1-1$.

\ \\Since $\deg_{x_j} f_i < d_j$ for $j=2,\ldots,n$ and $f_i(\beta_2,\ldots,\beta_n)=0$ for all $(\beta_2,\ldots,\beta_n)\in\prod_{i=2}^n V_i$, it follows by induction that $f_i = 0$.  This is true for each $i=0,\ldots,d_1 -1$.  So $f=0$.  $\Box$

\begin{lemma}
\label{Id}
Let $f_i (x_i)$ be a monic polynomial in $x_i$ of degree $d$ without multiple roots over $\bar{k}$, for $i=1,\ldots,n$.  Let $I\subset \bar{k} [x_1,\ldots,x_n]$ be the ideal generated by $f_i (x_i)$, $i=1,\ldots,n$.  Then $I_d$ is a $\bar{k}$-vector space with $f_i$, $i=1,\ldots,n$, as a basis.
\end{lemma}

\ \\{\bf Proof} Let $V_i$ be the set of roots of $f_i$ for $i=1,\ldots, n$. If $f\in I$ then $f$ vanishes at all points of $\prod_{i=1}^n V_i$.  Suppose $f\in I_d$.  Then $f$ as a polynomial in $\bar{k}[x_1,\ldots,x_n]$ of degree $d$ can be written in the form $f=\sum_{i=1}^n \alpha_i f_i (x_i) + g$ where $\alpha_i \in\bar{k}$ and $g\in \bar{k}[x_1,\ldots,x_n]$ where
$\deg_{x_i} g < d$ for $i=1,\ldots,n$.  Since $\sum_{i=1}^n \alpha_i f_i (x_i)\in I$ and $f\in I$, it follows $g\in I$ hence $g$ vanishes at all points of $\prod_{i=1}^n V_i$.  Since $\deg_{x_i} g < d$ for $i=1,\ldots,n$, it follows from Lemma~\ref{f=0} that $g=0$.  Therefore $f=\sum_{i=1}^n \alpha_i f_i (x_i)$ with $\alpha_i \in\bar{k}$.  $\Box$

\begin{lemma}
\label{IV}
Suppose $V,W\subset\bar{k}^n$ and $V=\lambda(W)$ where $\lambda\in Gl_n (\bar{k})$.  Then we have a $\bar{k}$-linear bijection
$I(V)\stackrel{\lambda^*}{\to} I(W)$ where $\lambda^*(f) = f\circ \lambda$ for $f\in I(V)$.
Moreover $\lambda^* (I(V)_i) = I(W)_i$ for all $i\ge 0$.
\end{lemma}

\ \\{\bf Proof} If $f\in I(V)$ then $f(V)=f(\lambda(W))=0$ hence $\lambda^* (f) \in I(W)$.  It is easily checked that $\lambda^*$ is $\bar{k}$-linear.  It is injective since $f\circ \lambda = 0$ implies $(f\circ\lambda)\circ \lambda^{-1} =f=0$.  It is surjective since for $g\in I(W)$, $g=f\circ\lambda$ where $f=g\circ\lambda^{-1} \in I(V)$.  Finally the map is degree preserving since $\lambda$ is linear.  Hence $\lambda^* (I(V)_i) = I(W)_i$ for all $i\ge 0$.  $\Box$

\begin{lemma}
\label{mixed}
Suppose $g(x)$ is univariate and monic of degree $d$ where $d$ is not divisible by the characteristic of the field $k$.   Suppose $\ell=\sum_{i=1}^n a_i x_i$ with $a_i\in\bar{k}$ for $i=1,\ldots,n$.
If $g\circ\ell = \sum_{i=1}^n \beta_i f_i (x_i)$ where $\beta_i\in\bar{k}$ for $i=1,\ldots,n$, then  $\ell = a_i x_i$ for some $i$.
\end{lemma}
\ \\{\bf Proof}  We will prove the lemma by induction on $n$.  The case $n=1$ is trivial.  Suppose inductively the lemma is true for $n-1$.   Let $\alpha_i$, $i=1,\ldots,d$, be the roots of $g$.  Then
$g\circ\ell = \prod_{i=1}^d (\ell-\alpha_i)$.  So we have
\[ \prod_{i=1}^d (\ell-\alpha_i) = \sum_{i=1}^n \beta_i f_i (x_i).\]
We see that for the left-hand-side of the equality,  the coefficient for $x_1^{d-1} x_2$ is $d a_1^{d-1} a_2$, which must be 0 by virtue of the right hand side.  Since $d$ is not divisible by the characteristic of $k$, we must have $a_1 a_2 = 0$.
Say $a_1=0$, then $\ell = a_2 x_2 +\ldots+ a_n x_n$ and $\beta_1 = 0$.  We are reduced to the case of $n-1$.  By induction we conclude that $\ell = a_i x_i$ for some $i\neq 1$.  Similarly if $a_2 = 0$ we conclude that $\ell = a_i x_i$ for some $i\neq 2$.  $\Box$

\ \\{\bf Proof of Theorem~\ref{main}}
Assume (2).  Then the space spanned by the $n$ polynomials $g_i\circ\mu$, $i=1,\ldots,n$ has dimension exactly $n$, since it contains $f_i (x_i)$, $i=1,\ldots,n$, which are clearly independent.  Moreover the space is also spanned by $f_i (x_i)$, $i=1,\ldots,n$.  Let $V_i$ be the set of roots of $f_i (x_i)$ for $i=1,\ldots,n$.  Let $V=\prod_{i=1}^n V_i$. Then by Lemma~\ref{Id}, $f_i (x_i)$, $i=1,\ldots,n$ generate $I(V)$ and they also form a basis  $I(V)_d$.
Let $\lambda = \mu^{-1}$.  Then $\alpha\in W$ if and only if $g_i (\alpha) = (g_i\circ\mu)(\lambda(\alpha))=0$ if and only if $f_i (\lambda (\alpha)) =0$, $i=1,\ldots,n$, if and only if $\lambda(\alpha) \in V$.  So $\lambda$ maps bijectively from $W$ to $V$.  So (2) implies (1).

Assume (1).  We claim that $g_1,\ldots,g_n$ are linearly independent.  Otherwise say $g_n$ is linearly dependent on $g_1,\ldots, g_{n-1}$.  Then $Z(\cG)=Z(\{g_1,\ldots,g_{n-1}\})$, which has dimension greater than 0, and we have a contradiction.

By Lemma~\ref{IV}, $I(W)_d$ has dimension equal to the dimension $I(V)_d$.  By Lemma~\ref{Id}, $I(V)_d$ has dimension $n$. Therefore $g_1,\ldots,g_n$ form a basis of $I(W)_d$.  Let $\mu=\lambda^{-1}$.  By Lemma~\ref{IV}, $g_i\circ\mu$, $i=1,\ldots,n$, form a basis of $I(V)_d$.  Let $V_i$ be the set of roots of $f_i (x_i)$ for $i=1,\ldots,n$.  Then
by Lemma~\ref{Id} $f_i (x_i)$, $i=1,\ldots,n$, form a basis of $I(V)_d$, which as argued above is also spanned by
$g_i\circ\mu$, $i=1,\ldots,n$.  So (1) implies (2).  $\Box$

\ \\{\bf Proof of Theorem~\ref{main-unique}}  Suppose $W\stackrel{\lambda}{\to} \prod_{i=1}^n V_i$ and $W\stackrel{\lambda'}{\to} \prod_{i=1}^n V'_i$
with $|V_i|=|V'_i| = d$ for $i=1,\ldots,n$ and $\lambda,\lambda'\in Gl_n (\bar{k})$.  Let $V=\prod_{i=1}^n V_i$ and
$V'=\prod_{i=1}^n V'_i$.  Then $V \stackrel{\rho}{\to} V'$ with $\rho=\lambda'\circ\lambda^{-1}\in Gl_n (\bar{k})$.
Let $f_i (x_i) = \prod_{\alpha\in V_i} (x_i -\alpha)$ and $g_i (x_i) = \prod_{\alpha\in V'_i} (x_i -\alpha)$ for $i=1,\ldots,n$.  Then $g_i \in I(V')_d$ and by Lemma~\ref{IV} $g_i\circ\rho\in I(V)_d$ for $i=1,\ldots,n$.  By Lemma~\ref{Id}, $g_i\circ\rho = \sum_{j=1}^n \beta_{ij} f_j (x_j)$ with $\beta_{ij}\in\bar{k}$, for $i,j=1,\ldots,n$.
Now the theorem follows easily from Lemma~\ref{mixed}.  $\Box$

\section{Last fall degree and Computation}
\label{compute}
Let $\mathcal{F}$ be a finite set of polynomials in  $k[x_1,\ldots,x_n]$.
For $i \in \Z_{\geq 0}$, we let $V_{\mathcal{F},i}$ be the smallest $k$-vector space such that
\begin{enumerate}
\item
$\{f \in \mathcal{F}: \deg(f) \leq i \} \subseteq V_{\mathcal{F},i}$;
\item
if $g \in V_{\mathcal{F},i}$ and if $h \in k[ x_1,\ldots,x_n]$ with $\deg(hg) \leq i$, then $hg \in V_{\mathcal{F},i}$.
\end{enumerate}

We write $f\equiv_i g \pmod{\mathcal{F}}$, for $f,g\in k[x_1,\ldots,x_n]$, if $f-g\in V_{\mathcal{F},i}$.

Put $R=k[x_1,\ldots,x_n]$ and let $R_{\le i}$ denote the $k$-vector space consisting of polynomials in $R$ of degree no greater than $i$ for $i \ge 0$.  The \emph{last fall degree} as defined in \cite{HKY} (see also \cite{HKYY}) is the largest $d$ such that $V_{\mathcal{F},d} \cap R_{\le d-1} \neq V_{\mathcal{F},d-1}$. We denote the last fall degree of $\mathcal{F}$ by $d_{\mathcal{F}}$.

As shown in \cite{HKY,HKYY} the last fall degree  is intrinsic to a polynomial
system, independent of the choice of a monomial order, always bounded by the degree of regularity,  and  invariant under linear change of variables and linear change of equations.

\begin{theorem}
\label{dG}
Suppose $W\subset\bar{k}^n$ and  $W\stackrel{\lambda}{\to} \prod_{i=1}^n V_i$ with $\lambda\in Gl_n (\bar{k})$, $V_i\subset\bar{k}$ and $| V_i | \le d$ for $i=1,\ldots, n$.  Suppose $\cF\subset \bar{k} [ x_1,\ldots,x_n]$, $Z(\cF) = W$ and the vector space spanned by $\lambda^*(\cF)$ contains $f_i$ for $i=1,\ldots,n$ where $f_i (x_i ) = \prod_{\alpha\in V_i} (x_i - \alpha)$. Then $d_{\cF}\le d$.
\end{theorem}
\ \\{\bf Proof}  The last fall degree is invariant under linear change of variables and linear change of equations (Proposition 2.6 (part v) of \cite{HKYY}).  Therefore we are reduced to the case that $W=\prod_{i=1}^n V_i$ and $\cF$ contains $f_i (x_i ) = \prod_{\alpha\in V_i} (x_i - \alpha)$ for $i=1,\ldots,n$.
We have
$f_i (x_i) = x_i^{d} - h_i (x_i)$ with $\deg h_i < d$.  So
$x_i^{d_i} \equiv_{d_i} h_i (x_i) \pmod{f_i}$.  Inductively we get
$x_i^j \equiv_j h_{ij} (x_i) \pmod{f_i}$ for some $h_{ij} (x_i)$ of degree less than $d_i$, for all $j \ge 0$.
From this it is easy to see that for all $f\in\bar{k} [ x_1\ldots x_n]$, we have
$f\equiv_{\max(d,\deg f)} \bar{f} \pmod{I}$ where $\bar{f}\in \bar{k}[x_1,\ldots, x_n]$ with $\deg_{x_i} \bar{f} < d_i$ for all $i$.
Finally $f\in I$ if and only if $\bar{f}\in I$, and since $\deg_{x_i} \bar{f} < d_i$ for all $i$, we conclude from
Lemma~\ref{f=0}
$\bar{f} \in I$ if and only if $\bar{f}=0$. It follows that $I_i = V_{\cF, i}$ for $i\ge d$.  Hence $d_{\cF} \le d$. $\Box$

\subsection{Last fall degree with field equations}
Let $k$ be a finite field with $q=p^m$ elements where $p$ is prime.
Suppose $\lambda\in Gl_n (k)$ and  $\cF\subset k[x_1,\ldots,x_n]$ is a finite set.
Let $\lambda^*(\cF) =\{ f\circ\lambda : f\in \cF\}$.
Let $\cE=\{x_i^q - x_i : i=1,\ldots,n\}$, the set of field equations for $k$.
Let $\vec{x} = (x_1,\ldots,x_n)$.  Suppose $\lambda=(a_{ij})$ with $1\le i,j\le n$ and $a_{ij}\in k$.  Let
$\ell_i (\vec{x}) =\sum_{j=1}^n a_{ij}x_j$ for $i=1,\ldots,n$, so that $\lambda (\vec{x}) = (\ell_i (\vec{x}))_{i=1}^n$.
Let $\lambda (\cE) = \{\ell_i (\vec{x}^q-\vec{x}): i=1,\ldots,n\}$.

\begin{lemma}
\label{dFq}
$d_{\lambda^*(\cF) \cup\cE} = d_{\cF\cup\cE}$.
\end{lemma}
\ \\{\bf Proof}  As noted before, the last fall degree is invariant under linear change of variables and linear change of equations (Proposition 2.6 (part v) of \cite{HKYY}).  So $d_{\lambda^*(\cF\cup\cE)}=d_{\cF\cup\cE}$.
Since
\[ (\vec{x}^q -\vec{x}) \circ\lambda = ((\ell_i (\vec{x}))^q - \ell_i (\vec{x}))_{i=1}^n =
(\ell_i (\vec{x}^q - \vec{x}))_{i=1}^n,\]
it follows that $\lambda^* (\cE) = \lambda (\cE)$.
Since $\cE$ and $\lambda(\cE)$ are related by a linear change of equations, so are
$\lambda^*(\cF)\cup\lambda (\cE)$ and $\lambda^*(\cF)\cup\cE$, so
$d_{\lambda^*(\cF)\cup\cE} = d_{\lambda^*(\cF)\cup\lambda(\cE)}$.  Therefore
\[  d_{\cF\cup\cE}= d_{\lambda^*(\cF\cup\cE)}= d_{\lambda^*(\cF)\cup\lambda^*(\cE)}= d_{\lambda^*(\cF)\cup\lambda (\cE)}=
d_{\lambda^*(\cF)\cup\cE}.\]  $\Box$

We observe that $Z(\cF \cup \cE) = Z_k (\cF)$ where $Z_k (\cF) = Z(\cF) \cap k^n$.  The degree of field equations for $k=\F_q$ can be reduced by introducing more variables.  More precisely $k=\F_q=Z(x^q-x)\simeq Z(\{ x^p-x_1, x_1^p-x_2,\ldots,x_{m-1}^p-x\})$ where $x\in k$ corresponds to $(x,x_1,\ldots,x_{m-1})$, with $x_1=x^p$,\ldots, $x_{m-1} = x_{m-2}^p$.  Now consider $nm$ variables $x_{ij}$, $i=1,\ldots,n$, $j=0,\ldots,m-1$, and identify $f(x_1,\ldots,x_n) \in\cF$ with $f(x_{10},\ldots, x_{n0})$.  Let $\cE'$ include the polynomials $x_{i0}^p-x_{i1}$,\ldots,$x_{i \ m-1}^p-x_{i0}$ for $i=1,\ldots,n$.  Then $Z(\cE') \simeq k^n$ where $(x_{ij})$ corresponds to $(x_{i0})$, and
$Z(\cF \cup \cE')\simeq Z_k (\cF)$.

Let  $\sigma$ denotes the Frobenius $p$-th power map: $\sigma (x)=x^p$ for $x\in\bar{k}$.
Suppose $\lambda\in Gl_n (k)$.  Let $\tilde{\lambda} \in Gl_{nm} (k)$ such that $\tilde{\lambda}$ is block-diagonally decomposed as $(\lambda_0,\lambda_1,\ldots,\lambda_{m-1})$ with $\lambda_j=\lambda^{\sigma^j}$ acting on the block
$(x_{1j},\ldots,x_{nj})$.  Let $\vec{x}_j = (x_{1j},\ldots,x_{nj})$, for $j=0,\ldots,m-1$.  Then
$x_{ij}\circ\tilde{\lambda} = ( \lambda_j (\vec{x}_j) )_i$, so $\vec{x}_j\circ\tilde{\lambda}  =\lambda_j (\vec{x}_j)$.

\begin{lemma}
\label{dFpm}
$d_{\cF\cup\cE'} = d_{\lambda^*{\cF}\cup\cE'}$
\end{lemma}
\ \\{\bf Proof}  Again, the last fall degree is invariant under linear change of variables and linear change of equations (Proposition 2.6 (part v) of \cite{HKYY}).  So $d_{\tilde{\lambda}^*(\cF\cup\cE')}=d_{\cF\cup\cE'}$.
Let $\vec{x}_i = (x_{1i},\ldots,x_{ni})$, for $i=0,\ldots,m-1$.
Let $\vec{y}_i = \lambda_i (\vec{x}_i )$ for $i=0,\ldots,m-1$.
Let $\cE''$ consists of $y_{ij}^p- y _{i \ j+1}$ for $i=1,\ldots, n$ and $j=0,\ldots, m-1$ (and where $j+1$ is taken
$\mod m$).
Then
\[ (\vec{x}_j^p - \vec{x}_{j+1} )\circ \tilde{\lambda} = (\lambda_j (\vec{x}_j ))^p - \lambda_{j+1}( \vec{x}_{j+1} )
=  \vec{y}_j^p - \vec{y}_{j+1}, \]
where the index $j+1$ is taken $\mod m$.
Hence $\tilde{\lambda}^* (\cE') = \cE''$.

Let $\mu=\lambda^{-1}$, and $\mu_i = \lambda_i^{-1}$ for $i=0,\ldots,m-1$.
Then for $j=0,\ldots,m-1$, $\vec{x}_j = \mu_j (\vec{y}_j)$ and
\[ \vec{x}_j^p - \vec{x}_{j+1}  = (\mu_j (\vec{y}_j ))^p - \mu_{j+1}( \vec{y}_{j+1} )
=  \mu_{j+1} (\vec{y}_j^p - \vec{y}_{j+1}), \]
where the index $j+1$ is taken $\mod m$.

We conclude that $\cE' = \tilde{\mu}(\cE'')$ where $\tilde{\mu}\in Gl_{mn} (k)$
is block-diagonally decomposed as $(\mu_0,\mu_1,\ldots,\mu_{m-1})$.

We have
$$d_{\cF\cup\cE'}=d_{\tilde{\lambda}^*(\cF\cup\cE')}=d_{\lambda^*(\cF)\cup\cE''}=d_{\lambda^*(\cF)\cup\tilde{\mu}(\cE'')}
= d_{\lambda^*(\cF)\cup\cE'}.$$  $\Box$

\subsection{Solving for $k$-rational points}
\begin{lemma}
\label{dp}
Let $k$ be a finite field with $|k| = q=p^m$.  Let $f\in k[x_0]$ with $d=\deg f$ and $\cF = \{ f, x_0^p-x_1,\ldots, x_{m-1}^p - x_0\}$.
Then $x_i \equiv_{pd} g_i $ with $g_i\in k[x_0]$ and  $\deg g_i < d$ for $i=1,\ldots,m-1$.  Suppose $g=gcd ( f, x_0^q-x_0 )\in k [ x_0]$.  Then $g\equiv_{dp} 0 \pmod{\cF}$.
\end{lemma}
\ \\{\bf Proof}  We have $x_0^p \equiv_{\max\{p,d\}} g_1 \pmod{f}$ for some $g_1\in k[x_0]$ with $\deg g_1 < d$.
We have $g_1^p \equiv_{pd} g_2 \pmod{f}$ for some $g_2\in k[x_0]$ with $\deg g_2 < d$.  Inductively we have
$g_i^p\equiv_{pd} g_{i+1} \pmod{f}$ with $g_{i+1}\in k[x_0]$ and $\deg g_{i+1} < d$, for $i=1,\ldots, m-2$.
In particular it follows that $x_0^q - x_0 \equiv h \pmod{f}$ where $h=g_m-x_0$, so $gcd (x_0^q-x_0,f)=gcd(f,h)$.

We have
$$x_1\equiv_p x_0^p \equiv_{\max\{p,d\}} g_1 \pmod{\cF}$$, and
inductively, $$x_{i+1} \equiv_p x_i^p \equiv_{pd} (g_i)^p\equiv_{pd} g_{i+1}\pmod{\cF},$$
for $i=1,\ldots, m-2$.
Finally, $x_0 \equiv_p x_{m-1}^p \equiv_{pd} g_{m-1}^p \equiv_{pd} g_m\pmod{\cF}$.
It follows that $h\equiv_{pd} 0 \pmod{\cF}$, consequently $gcd(h,f) \equiv_{pd} 0 \pmod{\cF}$.  Therefore
$gcd(x_0^q-x_0, f)\equiv_{pd} 0 \pmod{\cF}$.  $\Box$

\begin{theorem}
\label{dGE}
Let $k$ be a finite field with $|k| = q=p^m$.
Let
$\cG=\{g_1,\ldots,g_n\}$ where $g_i\in k[x_1,\ldots,x_n]$ with $\deg g_i \le d$ for $i=1,\ldots,n$.  Let $f_i\in k[x_i]$  for $i=1,\ldots,n$.  Let $\vec{g}=\left( \begin{array}{c} g_1\\. \\. \\g_n\end{array} \right)$ and  $\vec{f}=\left( \begin{array}{c} f_1\\. \\. \\f_n\end{array} \right)$.  Suppose $\vec{g}=\rho\circ\vec{f}\circ\lambda$ with
$\rho,\lambda\in Gl_n (k)$.
Then (1) $d_{\cG\cup\cE}\le d+q$ where $\cE=\{ x_i^q - x_i : i=1,\ldots, n\}$.    (2) $d_{\cG\cup\cE'}\le dp$, where $\cE'$ consists of $x_{ij}^p- x _{i \ j+1}$ for $i=1,\ldots, n$, $j=0,\ldots, m-1$ and $j+1$ is taken
$\mod m$, and $g(x_1,\ldots,x_n)\in k[x_1,\ldots,x_n]$ is identified with $g(x_{10},\ldots,x_{n0})\in k [x_{ij}: 1=1,\ldots,n, j=0,\ldots,m-1]$.
\end{theorem}
\ \\{\bf Proof}  Let $\cF=\{f_1,\ldots,f_n\}$ and $\mu=\lambda^{-1}$. Then $\vec{f}=\rho^{-1}\circ\vec{g}\circ\mu$.  We see that $\cF$ and $\mu^*\cG$ are related by $\rho^{-1}$. So by Lemma~\ref{dFq}, $d_{\cG\cup\cE}=d_{\mu^*{\cG}\cup\cE} = d_{\cF\cup\cE}$.  By Lemma~\ref{dFpm}, $d_{\cG\cup\cE'}=d_{\mu^*{\cG}\cup\cE'} = d_{\cF\cup\cE'}$.  So to prove (1) and (2) we are reduced to proving  $d_{\cF\cup\cE}\le d+q$ and  $d_{\cF\cup\cE'}\le dp$.

\ \\To prove $d_{\cF\cup\cE}\le d+q$,
let
$I'$ be the ideal generated by $\cF':=\cF\cup\cE$.
Let $h_i (x_i) = gcd (x_i^q - x_i, f_i (x_i))$ for $i=1,\ldots,n$.  Then $I'$ is generated by $\cH = \{h_i (x_i) : i=1,\ldots,n\}$.
Since $h_i = A (x_i^q-x_i) + B f_i$ with $A,B\in k[x_i]$ with $\deg A < \deg f_i$ and $\deg B < q$, $h_i (x_i) \equiv_{d_i + q} 0\pmod{\{x_i^q-x_i, f_i(x_i)\}}$.   So
$h_i (x_i) \equiv_{d_i+q} 0 \pmod{\cF'}$.
Let $c_i = \deg h_i$ for $i=1,\ldots,n$.
As with the proof of Theorem~\ref{dG},
it is easy to see that for all $f\in\bar{k} [ x_1\ldots x_n]$, we have
$f\equiv_{\max(c_i,\deg f)} \bar{f} \pmod{\cH}$ where $\bar{f}\in \bar{k}[x_1,\ldots, x_n]$ with $\deg_{x_i} \bar{f} < c_i$ for all $i$.
Finally $f\in I'$ if and only if $\bar{f}\in I'$, and since $\deg_{x_i} \bar{f} < d_i$ for all $i$, we conclude from
Lemma~\ref{f=0}
$\bar{f} \in I'$ if and only if $\bar{f}=0$. It follows that $I_i = V_{\cF', i}$ for $i\ge d+q$.  Hence $d_{\cF'} \le d+q$.

\ \\Next we prove $d_{\cF\cup\cE'}\le dp$ in a similar fashion by using Lemma~\ref{dp}. Put $\cF''=\cF\cup\cE'$ and let $J$ be the ideal generated by $\cF''$.  Let $\cH=\{ h_1,\ldots,h_n\}$ where $h_i = gcd (f_i, x_i^q-x_i)$ for $i=1,\ldots,n$.
By Lemma~\ref{dp} $\cH\subset V_{\cF'',dp}$.
Then
$J = \langle \cH\cup\cE' \rangle$ and $J\cap k [x_{i0}: i=1,\ldots, n] = \langle \cH \rangle$.
Moreover $Z(\cF'')$ maps bijectively to $Z(\cH)\subset k^n$ through the projection map sending $(x_{ij})_{i=1,\ldots,n; j=0,\ldots, m-1}$ to $(x_{i0})_{i=1,\ldots,n}$.
By Lemma~\ref{dp}, $x_{ij}\equiv_{pd} g_{ij}\pmod{\cF''}$ with $g_{ij}\in k[x_{i0}]$ and $\deg g_{ij} < d$.  We have
$x_{ij}^2 \equiv_{pd} g_{ij}^2\pmod{\cF''}$, and $g_{ij}^2 \equiv_{2d} g_{ij,2} \pmod{h_i}$ with $g_{ij,2}\in k[x_{i0}]$ and $\deg g_{ij,2} < d$.  Inductively we see that $x_{ij}^e \equiv_{pd} g_{ij,e} \pmod{\cF''}$ with $g_{ij,e}\in k[x_{i0}]$ and $\deg g_{ij,e} < d$, for $e < p$.
From this it is easy to see that for all $f\in k[x_{ij}: i=1,\ldots,n; j=0,\ldots,m-1]$,
$f\equiv_{\max(\deg f, pd)} \bar{f}$ where $\bar{f} \in k[ x_{i0}: i=1,\ldots,n]$ and $\deg_{x_{i0}} (\bar{f}) < \deg h_i$ for all $i$.  So $f\in J$ if and only if $\bar{f}\in J\cap k [x_{i0}: i=1,\ldots, n]=\langle \cH\rangle$.  We conclude from
Lemma~\ref{f=0}
$\bar{f} \in \langle \cH \rangle$ if and only if $\bar{f}=0$. It follows that $J_i = V_{\cF'', i}$ for $i\ge pd$.  Hence $d_{\cF'} \le pd$.
$\Box$


\end{document}